\documentclass[12pt]{article}
\usepackage{amsfonts}
\usepackage{}
\usepackage{mathrsfs,graphicx,amsmath,amssymb,amsfonts}

\usepackage{appendix}
\usepackage{caption}
\usepackage{mathrsfs,graphicx,amsmath,amssymb,amsfonts,amsthm}
\usepackage{cite}
\usepackage{longtable}
\usepackage{supertabular}
\evensidemargin=\oddsidemargin\topmargin=-1.5cm

\theoremstyle{definition}

\begin{document}
\title{The least eigenvalue of the complements of graphs with given connectivity \thanks{This work was supported by NSFC (No. ????).}}
\author{Huan Qiu$^1$, Keng Li$^2$, Jianfeng Wang$^1$\footnote{Corresponding author. Email: jfwang4@yahoo.com.cn.}, Guoping Wang$^3$\\
{\small 1.School of Mathematics and Statistics, Shandong University of Technology,}\\
{\small Zibo, Shandong 255000,P.R.China.}\\
{\small 2.College of Digital Economy, Hubei University of Automotive Technology,}\\
{\small Shiyan, Hubei 442000, P.R.China }\\
{\small 3.School of Mathematical Sciences, Xinjiang Normal University,}\\
{\small Urumqi, Xinjiang 830054, P.R.China. }\\
{\small 4.School of Mathematics and Statistics, Shandong University of Technology,}\\
{\small Zibo, Shandong 255000, P.R.China.}\\}

\date{}
\maketitle {\bf Abstract.}
The least eigenvalue of a graph $G$ is the least eigenvalue of the adjacency matrix of $G$.
In this paper, we determine the graphs which attain the minimum least eigenvalue
among all complements of connected simple graphs with given connectivity.

{\flushleft{\bf Key words:}} The least eigenvalue; Complements of graphs; Connectivity.\\
{\flushleft{\bf MR(2020) Subject Classification:}}  05C40, 05C50\\

\section{Introduction}
~~~~~Let $G$ be a simple graph with the vertex set $V(G)=\{v_1,v_2,\ldots,v_n\}$
and the edge set $E(G)$.
The adjacency matrix of $G$ is denoted by $A(G)=(a_{ij})_{n\times n}$,
where $a_{ij}=1$ if $v_iv_j\in E(G)$, and $a_{ij}=0$ otherwise.
Since $A(G)$ is a non-negative real symmetric matrix,
its eigenvalues can be ranged as $\lambda_1(G) \geq \lambda_2(G) \geq \cdots \geq \lambda_n(G)$,
where $\lambda_1(G)$ and $\lambda_n(G)$ are called spectral radius and the least eigenvalue of $G$, respectively.	

There is a lot of research on the least eigenvalue of a graph.
One of the goals is to use eigenvalue or eigenvector to investigate the graph structure.
Bell, Cvetkovi\'{c}, Rowlinson and Simi\'{c}\cite{F.B.1} characterized the graphs with minimum least eigenvalue.
Fan, Wang and Gao\cite{Y.Y.Y} determined the unique graph with minimum least eigenvalue among all unicyclic graphs with fixed order.
Further results on the least eigenvalue were focused on graphs with some invariants being fixed,
like connectivity by Ye, Fan and Liang\cite{M.L},
domination number, maximum vertex degree, cut vertices and cut edge by Zhu\cite{B.Z.O, B.Z,Y.F, B.Z.T},
independence number or cover number by Tan and Fan\cite{Y.T},
or some specified classes of graphs,
like unicyclic graphs with $k$ pendant vertices and  fixed diameter by Liu, Zhai and Shu\cite{R.S,M.S} and unicyclic graphs with given independence number by Hou and Qu\cite{H.Q}.

The complement of a graph $G$ is denoted by $G^c=(V(G^c),E(G^c))$,
where $V(G^c)=V(G)$ and $E(G^c)=\{uv:u,v\in V(G),uv\notin E(G)\}$.
A graph and its complement could be completely distinct,
say, the complement of a tree is no longer a tree, and so the least eigenvalues of a graph and its complement could be completely distinct.
Therefore, it is meaning for us to study the least eigenvalues of the complements of graphs.	
Up to now there is little research about the least eigenvalues of the complements of graphs.
Fan, Zhang and Wang \cite{Y.Fan} characterized the connected graph with the minimal least eigenvalue among all complements of trees.
Jiang, Yu, Sun and Ruan \cite{G.Jiang} achieved the graph with the minimum least eigenvalue among all graphs whose complements have only two pendent vertices.
	
The {\it connectivity} of the graph $G$ is the minimum number of vertices whose deletion yields the resulting graph disconnected. In this paper, we investigate the least eigenvalue of the complements of graphs with given connectivity.
It can be seen from the proof of our main theorem that our approach is quite different from those previously exploited in problems regarding the least eigenvalue of graphs. This enables us to determine the graph which attains the minimum least eigenvalue
among all complements of connected simple graphs with given connectivity.
	
\section{Preliminary}

Suppose that $G$ is a simple graph
with the vertex set $V(G)=\{v_1,v_2,\ldots, v_n\}$.
Let $x=(x_1,x_2,\cdots,x_n)^T$,
where $x_i$ corresponds to $v_i$, i.e., $x(v_i)=x_i$ for $i=1,2,\cdots,n$.
Then
	\begin{equation}
		x^TA(G)x=\sum_{v_iv_j\in E(G)}2x_ix_j.
	\end{equation}

The set of neighbours of $v$ in $G$ is denoted by $N_G(v)$.
Suppose that $x$ is an eigenvector of $A(G)$ corresponding to the eigenvalue $\lambda$.
Then for $v_i\in V(G)$, we have
	\begin{equation}
			\lambda x_i=\sum_{v_j\in N_G(v_i)}x_j\quad for\  i=1,2,\cdots,n.
	\end{equation}\vskip 3mm

{\noindent \bf Lemma 2.1.} (Rayleigh's inequalities)
{\it Let $G$ be a graph with spectral radius $\lambda_1(G)$ and the least eigenvalue $\lambda_n(G)$ of $A(G)$,
and $x=(x_1,x_2,\cdots,x_n)^T$ be a unit vector.
Then we have
\begin{equation}
	\lambda_n(G)\leq  x^TA(G)x\leq \lambda_1(G).
	\nonumber
\end{equation}
Moreover, the first equality holds if and only if $x$ is a unit eigenvector of $A(G)$ with respect to $\lambda_n(G)$
and the second equality holds if and only if $x$ is a unit eigenvector of $A(G)$ with respect to $\lambda_1(G)$.}\vskip 3mm

{\noindent \bf Lemma 2.2.} \cite{H.M} {\it Let $\Delta(G)$ be the maximum degree of a graph $G$.
Then $\Delta(G)\geq \lambda_1(G)\geq \sqrt{\Delta(G)}$.} \vskip 3mm

{\noindent \bf Lemma 2.3.} \cite{B.Z} {\it  Let $G^*$ be a connected graph with two non-adjacent vertices $u, v$ and let $G$ be the graph obtained from $G^*$ by adding the edge $uv$. Assume that $x$ and $y$ are the unit least vectors of $G$ and $G^*$, respectively.
Then

{\noindent \rm(i.)} $\lambda_n(G^*)\leq \lambda_n(G)$ if $x_u=0$ or $x_v=0$, and the equality holds if and only if x is a least vector of
$G^*$ and $x_u=x_v=0$.

{\noindent \rm(ii.)} $\lambda_n(G)\leq \lambda_n(G^*)$ if $y_u=0$ or $y_v=0$, and the equality holds if and only if y is a least vector of
$G$ and $y_u=y_v=0$.

{\noindent \rm(iii.)} $\lambda_n(G)< \lambda_n(G^*)$ if $y_uy_v< 0$. }

A {\it matching} in a graph is a set of pairwise nonadjacent edges.
A matching with $\alpha$ pairwise nonadjacent edges is denoted by $\alpha$-matching.\vskip 3mm

{\noindent \bf Lemma 2.4.} (Hall's theorem)
{\it A bipartite graph $G=G[X,Y]$ has a matching which covers every vertex in $X$ if and only if $|N_G(S)|\geq |S|$ for all $S\subseteq X$,
where $N_G(S)=\cup_{v\in S} N_G(v)$.}\vskip 3mm

From this lemma we easily see that the below result is true.\vskip 2mm

{\noindent \bf Corollary 2.5.} {\it Let $G$ be a graph with the vertex set $V(G)$.
Suppose $U$ and $W$ are two subsets of $V(G)$ satisfying $U\cap W=\emptyset$.
Then $G$ has a matching between $U$ and $W$ which covers every vertex in $U$ if $|N_G(S)\cap W|\geq |S|$ for all $S\subseteq U$.}\vskip 3mm

If neither its origin nor its terminus is covered by a matching $M$, the path is called an {\it M-augmenting path}.\vskip 3mm

{\noindent \bf Lemma 2.6.} (Berge's theorem)
{\it A matching $M$ in a graph $G$ is a maximum matching if and only if $G$ contains no $M$-augmenting path.}\vskip 3mm
	
{\noindent \bf Lemma 2.7.} \cite{Feng} {\it Suppose $G$ and $G^c$ are both connected graphs on $n$ vertices.
If $x$ is a least eigenvector of $G^c$
then $x$ contains at least two positive entries and negative entries.}

\section{Main results}
~~~Let $\mathcal{G}_{n,\kappa}$ denote the set of the connected simple graphs
on $n$ vertices with the connectivity $\kappa$,
and $\mathcal{G}_{n,\kappa}^c$ be the set of the complements of the graphs in $\mathcal{G}_{n,\kappa}$.
In this paper we usually consider $\mathbf{G}$ to be the graph of $\mathcal{G}_{n,\kappa}$
so that $\lambda_n(\mathbf{G}^c)$ is as small as possible in $\mathcal{G}_{n,\kappa}^c$.
If $\kappa=n-1$ then $\mathbf{G}$ is isomorphic to the complete graph $K_n$ of order $n$,
and so in what follows we assume $\kappa\leq n-2$.
Next we will use four claims to characterize $\mathbf{G}$.

In this paper we usually let $\partial(\mathbf{G})$ be a minimum vertex-cut of $\mathbf{G}$
and $x=(x_1,\cdots,x_n)^T$ be a unit eigenvector of $A(\mathbf{G}^c)$ with respect to $\lambda_n(\mathbf{G}^c)$. \vskip 3mm

Let $G$ be a graph on $n$ vertices,
$J_n$ be the matrix of order $n$ whose all entries are $1$
and $I_n$ be the identity matrix of order $n$.
Then we have
\begin{equation}
	A(G^c)=J_n-I_n-A(G).
	\nonumber
\end{equation}

{\noindent \bf Claim 3.1.} {\it $\mathbf{G}-\partial(\mathbf{G})$ contains exactly two components $\mathbf{G}_1$ and $\mathbf{G}_2$.}

{\noindent \bf Proof.}
Suppose on the contrary that $\mathbf{G}-\partial(\mathbf{G})$ contains three components $\mathbf{G}_1$, $\mathbf{G}_2$ and $\mathbf{G}_3$.
Obviously, the positive and negative values of vertices in $\mathbf{G}_1$ and $\mathbf{G}_2$ are not the same.
Then there must be two vertices $u\in V(\mathbf{G}_1)$ and $v\in V(\mathbf{G}_2)$ such that $x(u)x(v)\geq 0$.
Let $H_1=\mathbf{G}+uv$.
Note that $\partial(\mathbf{G})$ is also $\kappa$-vertex cut of $H_1$, and so $H_1\in \mathcal{G}_{n,\kappa}$.
From the equation $(1)$, we have
$x^TA(\mathbf{G})x=\sum_{v_iv_j\in E(\mathbf{G})}2x_ix_j \leq \sum_{v_iv_j\in E(H_1)}2x_ix_j=x^TA(H_1)x$.

Thus, by Lemma 2.1, we have
\begin{equation}
	\begin{split}
		\lambda _n(\mathbf{G}^c) &= x^TA ( \mathbf{G}^c  ) x\\
		&=x^T ( J_n-I_n  ) x-x^TA ( \mathbf{G}  ) x\\
		&\geq x^T ( J_n-I_n  ) x-x^TA (H_1) x\\
		&= x^TA(H_1^c)x\\
		&\geq \lambda _n(H_1^c).
	\end{split}
	\nonumber
\end{equation}
Lemma $2.3$ shows the equality holds if and only if $x$ is a least vector of $H_1^c$ and $x(u)=x(v)=0$.
If $\lambda _n(\mathbf{G}^c) =\lambda _n(H_1^c)$ then we consider $H_1$ as $\mathbf{G}$,
and otherwise $\lambda_n(\mathbf{G}^c)> \lambda_n(H_1^c)$,
which contradicts the choice of $\mathbf{G}^c$,
and so Claim 3.1 is true. $\Box$ \vskip 3mm

Suppose $G$ is a graph with the vertex set $V(G)$.
If $S$ is a non-empty subset of $V(G)$ then
we denote by $G[S]$ the subgraph of $G$ induced by $S$.
Set $\partial(\mathbf{G})^+=\{v\in \partial(\mathbf{G}): x(v)\geq 0\}$, $\partial(\mathbf{G})^-=\{v\in \partial(\mathbf{G}):x(v)<0\}$,
$V_i^+=\{v\in V(\mathbf{G}_i):x(v)\geq 0\}$ and
$V_i^-=\{v\in V(\mathbf{G}_i):x(v)<0\}$, where $i=1,2$. \vskip 3mm

{\noindent \bf Claim 3.2.}
{\it $G[V_1^+\cup \partial(\mathbf{G})^+]$,
	$G[V_2^+\cup \partial(\mathbf{G})^+]$,
	$G[V_1^-\cup \partial(\mathbf{G})^-]$ and
	$G[V_2^-\cup \partial(\mathbf{G})^-]$ are all complete graphs.}

{\noindent \bf Proof.}
Suppose on the contrary that the vertices $u$ and $v$ are not adjacent in $G[V_1^+\cup \partial(\mathbf{G})^+]$.
Let $H_2=\mathbf{G}+uv$.
Note that $\partial(\mathbf{G})$ is also $\kappa$-vertex cut of $H_2$,
and so $H_2\in \mathcal{G}_{n,\kappa}$.
From the equation $(1)$,
$x^TA(\mathbf{G})x=\sum_{v_iv_j\in E(\mathbf{G})}2x_ix_j \leq \sum_{v_iv_j\in E(H_2)}2x_ix_j=x^TA(H_2)x$.
As in the proof of Claim 3.1,
we can verify $\lambda_n(\mathbf{G}^c)\geq \lambda_n(H_2^c)$.
Lemma $2.3$ shows the equality holds if and only if $x$ is a least vector of $\mathbf{G}^c$ and $x(u)=x(v)=0$.
If $x(u)=x(v)=0$ then we consider $H_2$ as $\mathbf{G}$,
and otherwise $\lambda_n(\mathbf{G}^c)> \lambda_n(H_2^c)$,
which contradicts the choice of $\mathbf{G}^c$,
and so $u$ and $v$ are adjacent in $G[V_1^+\cup \partial(\mathbf{G})^+]$.
Therefore, $G[V_1^+\cup \partial(\mathbf{G})^+]$ is complete graph.

Similarly, $G[V_2^+\cup \partial(\mathbf{G})^+]$,
$G[V_1^-\cup \partial(\mathbf{G})^-]$ and
$G[V_2^-\cup \partial(\mathbf{G})^-]$ are also complete graphs. $\Box$\vskip 3mm

Set $V^+=\{v\in V(\mathbf{G}): x(v)\geq 0\}$,
and $V^-=\{v\in V(\mathbf{G}):x(v)< 0\}$.
Write $|V^+|=n_1$ and $|V^-|=n_2$. Then $n_1+n_2=n$.
For convenience we assume without loss of generality that $n_1\geq n_2$,
in which case we can distinguish three cases as follows:\\
$$n_2\geq \kappa, ~~n_1\geq \kappa > n_2, ~~\kappa>n_1.$$
\vskip 3mm

Let $U$ be the subset of $V^+$ containing such vertices of $V^+$ that connect at least one vertex in $V^-$.
Let $W$ be the subset of $V^-$ containing such vertices of $V^-$ that connect at least one vertex in $V^+$.\vskip 3mm

{\noindent \bf Claim 3.3.} {\it  If $n_2\geq \kappa$, then there is a $\kappa$-matching between $V^+$ and $V^-$.}

{\noindent \bf Proof.}
Let $M$ be a maximum matching between $U$ and $W$.
Since $n_1\geq n_2$, $n_1\geq \kappa$.
If $|U|< \kappa$, then $\mathbf{G} [V(\mathbf{G})\backslash U]$ is not connected, and so $|U|\geq \kappa$.
Similarly,  $|W|\geq \kappa$.
Denote by $V(M)$ the set of the vertices which are covered by $M$.
If $U\backslash V(M)=\emptyset$ or $W\backslash V(M)=\emptyset$,
then there is an $\kappa$-matching between $U$ and $W$.

So we assume $U\backslash V(M)\neq \emptyset$ and $W\backslash V(M)\neq \emptyset$.
In this case there is no edge between $U\backslash V(M)$ and $W\backslash V(M)$
otherwise the edge is not matched under $M$,
which contradicts the maximality of $M$.

Let $S_1$ be the set of the vertices of $U\cap V(M)$ which are adjacent to some vertices of $W\backslash V(M)$.
Let $T_2$ be the set of the vertices of $W\cap V(M)$ which are adjacent to some vertices of $U\backslash V(M)$.
Let $T_1$ be the set of the vertices of $W\cap V(M)$ which are matching with the vertices of $S_1$.
Let $S_2$ be the set of the vertices of $U\cap V(M)$ which are matching with the vertices of $T_2$.
Then there is no edge between $T_1$ and $S_2$
otherwise we will get an $M$-augmenting path $P={w_1s_1t_1s_2t_2u_1}$,
where $w_1\in W\backslash V(M)$, $s_1\in S_1$, $t_1\in T_1$, $s_2\in S_2$, $t_2\in T_2$, $u_1\in U\backslash V(M)$ and $s_1t_1$, $s_2t_2\in M$.
This contradicts Lemma 2.6.

Clearly, $T_1\cap T_2=\emptyset$ and $S_1\cap S_2=\emptyset$.
Let $T_3$ be the set of the vertices of $(W\cap V(M))\backslash (T_1\cup T_2)$ such that for each vertex $t_3\in T_3$,
there is a path from $t_3$ to some vertex of $T_2$ whose edges are alternately in $\widetilde{E}\backslash M$ and $M$,
where $\widetilde{E}$ is the set of the edges between $(U\cap V(M))\backslash S_1$ and $(W\cap V(M))\backslash (T_1\cup T_2)$.
Let $S_3$ be the set of the vertices of $(U\cap V(M))\backslash (S_1\cup S_2)$ which are matching with the vertices of $T_3$.
Then there is no edge between $T_1$ and $S_3$ otherwise we will get an $M$-augmenting path.
This contradicts Lemma 2.6.

Let $S_4$ be the set of the vertices of $(U\cap V(M))\backslash (S_1\cup S_2)$ such that for each vertex $s_4\in S_4$,
there is a path from $s_4$ to some vertex of $S_1$ whose edges are alternately in $\widetilde{E}^*\backslash M$ and $M$,
where $\widetilde{E}^*$ are the set of the edges between $(W\cap V(M))\backslash T_2$ and $(U\cap V(M))\backslash (S_1\cup S_2)$.
Let $T_4$ be the set of the vertices of $(W\backslash V(M))\backslash (T_1\cup T_2)$ which are matching with the vertices of $S_4$.
Just as the above argument we can verify that there is no edge between $T_4$ and $S_2\cup S_3$.

$S_3\cap S_4=\emptyset$ and $T_3\cap T_4=\emptyset$ otherwise we will get an $M$-augmenting path,
which contradicts Lemma 2.6.
Let $S_5=(U\cap V(M))\backslash \bigcup _{i=1}^4S_i$ and $T_5=(W\cap V(M))\backslash \bigcup _{i=1}^4T_i$.
Clearly, the vertices of $S_5$ are matching with the vertices of $T_5$.
$T_i$ and $S_i$ ($1\leq i\leq 5$) are shown in Figure 1.
\begin{center}
\begin{figure}[htbp]
\centering
\includegraphics[height=30mm]{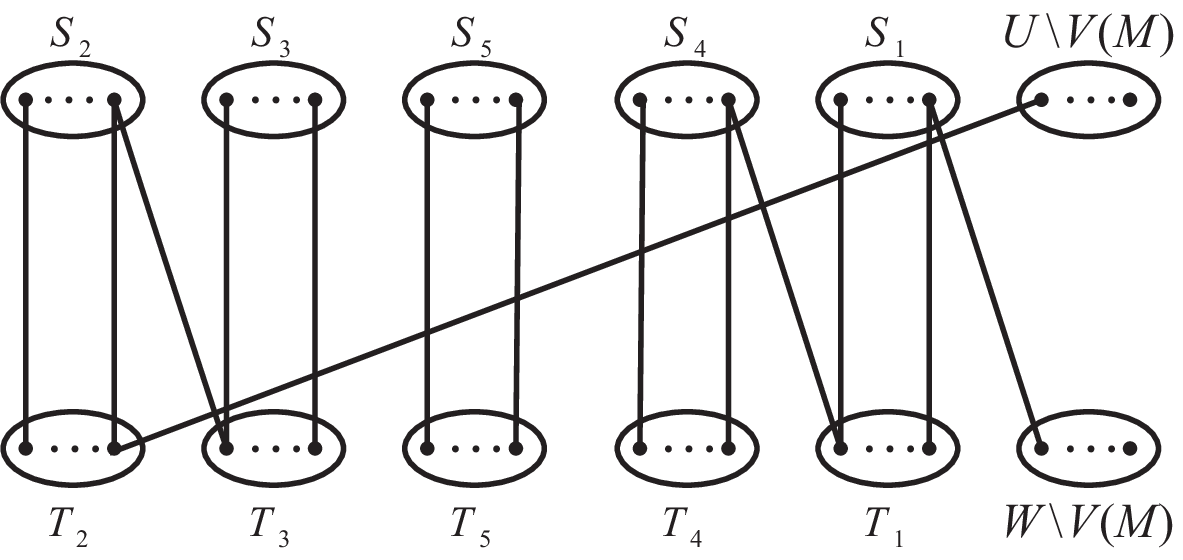}
\\Fig. 1.  $T_i$ and $S_i$ ($1\leq i\leq 5$)
\end{figure}
\end{center}

There is no edge between $T_5$ and $S_2\cup S_3$
otherwise we can find one path from the vertex $t_5$ of $T_5$ to some vertex of $T_2$
whose edges are alternatively in $\widetilde{E}\backslash M$ and $M$, which implies $t_5\in T_3$, a contradiction.
Similarly, we can verify that there is no edge between $S_5$ and $T_1\cup T_4$.

In summary, there is no edge between $T_1$ and $S_2\cup S_3\cup S_5\cup (U\cap V(M))$, between $T_4$ and $S_2\cup S_3\cup S_5\cup (U\cap V(M))$, and between $T_5$ and $S_2\cup S_3\cup (U\cap V(M))$.
Suppose that $M$ is an $\ell$-matching.
Then it is clear that $|S_1\cup T_2\cup T_3\cup S_4\cup S_5|=\ell$.
From the above argument, we easily see that $S_1\cup T_2\cup T_3\cup S_4\cup S_5$ is a vertex cut of $\mathbf{G}$,
and so $\ell\geq \kappa$.   $\Box$\vskip 3mm

Let $n_1\geq n_2\geq \kappa$ and $K_{n_1}$ and $K_{n_2}$ be two disjoint complete graphs.
Then we denote by $\mathbf{B}_1(n_1,n_2; \kappa)$ the  graph obtained from $K_{n_1}$ and $K_{n_2}$
by connecting $\kappa$ edges between $V(K_{n_1})$ and $V(K_{n_2})$ so that they become an $\kappa$-matching.
Clearly, $\mathbf{B}_1(n_1,n_2; \kappa)\in \mathcal{G}_{n,\kappa}$.\vskip 3mm


{\noindent \bf Lemma 3.1.} {\it When $n_2\geq \kappa$,
$\lambda_n(\mathbf{G}^c)\geq \lambda_n((\mathbf{B}_1^c(n_1,n_2;\kappa)))$.}

{\noindent \bf Proof.} From the claim 3.3 we know that there is an $\kappa$-matching $M$ between $V^+$ and $V^-$.
Connect respectively all pairs between $V_1^+$ and $V_2^+$ and between $V_1^-$ and $V_2^-$,
and delete all edges between $V^+$ and $V^-$ but those edges matched under $M$.
Then we obtain the resulting graph which is isomorphic to $\mathbf{B}_1(n_1,n_2;\kappa)$.

From the equation $(1)$, $x^TA(\mathbf{G})x=\sum_{v_iv_j\in E(\mathbf{G})}2x_ix_j \leq \sum_{v_iv_j\in E(\mathbf{B}_1(n_1,n_2;\kappa))}2x_ix_j=x^TA(\mathbf{B}_1(n_1,n_2;\kappa))x$.
As in the proof of Claim $3.1$, we can verify that $\lambda_n(\mathbf{G}^c)\geq \lambda_n((\mathbf{B}_1^c(n_1,n_2;\kappa)))$.  $\Box$\vskip 3mm

{\noindent \bf Claim 3.4.} {\it  If $n_1\geq \kappa>n_2$,
then there is a $n_2$-matching between $n_2$ vertices of $U$ and all vertices of $V^-$
and the other $\kappa-n_2$ vertices of $U$ are adjacent to each vertex of $V^-$.}

{\noindent \bf Proof.} If $V^-\backslash W\neq \emptyset$,
then $W$ is clearly a vertex cut of $\mathbf{G}$.
Whereas $|W|< \kappa$, this contradiction shows $W= V^-$.
If $|U|< \kappa$ then $U$ is clearly a vertex cut of $\mathbf{G}$.
This contradiction shows $|U|\geq \kappa$.

Set $U=\{u_1,u_2,\cdots,u_{|U|}\}$, where $x(u_i)\leq x(u_{i+1})$ ($1\leq i\leq |U|-1$),
and $R_1=\{u_1,u_2,\cdots,u_{\kappa-n_2}\}$ is subset of $U$.
If there is some $Q^*\subseteq V^-$ such that $|N_\mathbf{G}(Q^*)\cap (U\backslash R_1)|< |Q^*|$,
then $R_1\cup (V^-\backslash Q^*) \cup (N_\mathbf{G}(Q^*)\cap (U\backslash R_1))$ is a vertex cut of $\mathbf{G}$.
Whereas $|R_1\cup (V^-\backslash Q^*) \cup (N_\mathbf{G}(Q^*)\cap (U\backslash R_1))|< \kappa$,
this contradiction shows that for any $Q\subseteq V^-$, $|N_\mathbf{G}(Q)\cap (U\backslash R_1)|\geq |Q|$.
By Corollary 2.5, there exists a $n_2$-matching $M_1$ between $V^-$ and $U\backslash R_1$.

If there is some vertex $v_1\in V^-$ such that $|N_\mathbf{G}(v_1)\cap U|<\kappa-n_2+1$,
then $(V^-\backslash \{v_1\}) \cup (N_\mathbf{G}(v_1)\cap U)$ is a vertex cut of $\mathbf{G}$.
Whereas $|(V^-\backslash \{v_1\}) \cup (N_\mathbf{G}(v_1)\cap U)|< \kappa$,
this contradiction shows that for any $v\in V^-$, $|N_\mathbf{G}(v)\cap U|\geq \kappa-n_2+1$,
and so $v$ connects at least $\kappa-n_2$ vertices of $U$ except the vertex matching with $v$.
Now we will show that $v$ connects all vertices of $R_1$.

Suppose for a contradiction that some vertex $u_{\ell_1}$ of $R_1$ is not adjacent to $v$, where $1\leq \ell_1\leq \kappa-n_2$.
From the above argument, we easily see that $|N_\mathbf{G}(v)\cap U|\geq \kappa-n_2+1$.
Then there must be the vertex $u_{\ell_2}$ of $U\backslash R_1$ satisfying $vu_{\ell_2}\not\in M_1$,
where $\kappa-n_2\leq \ell_2\leq |U|$.
We delete the edge $vu_{\ell_2}$ and add $vu_{\ell_1}$.
Note $x(u_{\ell_1})\leq x(u_{\ell_2})$.
From the equation $(1)$,
$x^TA(\mathbf{G})x=\sum_{v_iv_j\in E(\mathbf{G})}2x_ix_j \leq \sum_{v_iv_j\in E(\mathbf{G}-vu_{\ell_2}+vu_{\ell_1})}2x_ix_j$.
As in the proof of Claim $3.1$,
we can verify $\lambda_n(\mathbf{G}^c)\geq \lambda_n((\mathbf{G}-vu_{\ell_2}+vu_{\ell_1})^c)$.
This contradiction shows $v$ connects all vertices of $R_1$. $\Box$\vskip 3mm

Let $n_1\geq \kappa> n_2$. Suppose $K_{n_1}$ and $K_{n_2}$ are disjoint.
Then we denote by $\mathbf{B}_2(n_1,n_2; \kappa)$ the  graph obtained from $K_{n_1}$ and $K_{n_2}$
by connecting $n_2$ edges between $V(K_{n_1})$ and $V(K_{n_2})$ so that they become an $n_2$-matching $M_1$
and connecting each vertex of $V(K_{n_2})$ and all $\kappa-n_2$ vertices of $V(K_{n_1})$ which are not covered by $M_1$.
Clearly, $\mathbf{B}_2(n_1,n_2; \kappa)\in \mathcal{G}_{n,\kappa}$.\vskip 3mm

{\noindent \bf Lemma 3.2.} {\it When $n_1\geq \kappa> n_2$,
$\lambda_n(\mathbf{G}^c)\geq \lambda_n((\mathbf{B}_2^c(n_1,n_2;\kappa)))$.}

{\noindent \bf Proof.}
From the claim 3.4 we easily observe that after connecting respectively all pairs of vertices between $V_1^+$ and $V_2^+$ and between $V_1^-$ and $V_2^-$
and deleting some edges between $V^+$ and $V^-$ we can
obtain the resulting graph which is isomorphic to $\mathbf{B}_2(n_1,n_2;\kappa)$.

From the equation $(1)$,
$x^TA(\mathbf{G})x=\sum_{v_iv_j\in E(\mathbf{G})}2x_ix_j \leq \sum_{v_iv_j\in E(\mathbf{B}_2(n_1,n_2;\kappa))}2x_ix_j=x^TA(\mathbf{B}_2(n_1,n_2;\kappa))x$.
As in the proof of Claim $3.1$,
we can verify that $\lambda_n(\mathbf{G}^c)\geq \lambda_n(\mathbf{B}_2^c(n_1,n_2;\kappa))$.
$\Box$\vskip 3mm

Let $\kappa>n_1$.
Suppose $K_{n_1}$ and $K_{n_2}$ are disjoint.
Set $S$ to be a subset of $V(K_{n_1})$ such that $|S|=n_1-n_2$.
Then we denote by $\mathbf{B}_3(n_1,n_2;\kappa)$ the graph obtained from $K_{n_1}$ and $K_{n_2}$ by connecting
each vertex of $S$ with each vertex of $V(K_{n_2})$
and connecting each vertex of $V(K_{n_1})\backslash S$ with $\kappa-n_1+1$ vertices of $V(K_{n_2})$
so that for any two vertices $s_1$ and $s_2$ of $V(K_{n_1})\backslash S$,
$N_{\mathbf{B}_3(n_1,n_2;\kappa)}(s_1)\cap V(K_{n_2})\neq N_{\mathbf{B}_3(n_1,n_2;\kappa)}(s_2)\cap V(K_{n_2})$,
and connecting each vertex of $V(K_{n_2})$ with $\kappa-n_1+1$ vertices of $V(K_{n_1})\backslash S$
so that for any two vertices $t_1$ and $t_2$ of $V(K_{n_2})$,
$N_{\mathbf{B}_3(n_1,n_2;\kappa)}(t_1)\cap V(K_{n_1})\neq N_{\mathbf{B}_3(n_1,n_2;\kappa)}(t_2)\cap V(K_{n_1})$.
Clearly, $\mathbf{B}_3(n_1,n_2; \kappa)\in \mathcal{G}_{n,\kappa}$.
\vskip 3mm

{\noindent \bf Lemma 3.3.} {\it When $\kappa>n_1$,
	$\lambda_n(\mathbf{G}^c)\geq \lambda_n((\mathbf{B}_3^c(n_1,n_2;\kappa)))=\kappa+1-n$.}

{\noindent \bf Proof.}
If $V^+\backslash U\neq \emptyset$, then $U$ is clearly a vertex cut of $\mathbf{G}$.
Whereas $|U|<\kappa$, this contradiction shows $U=V^+$.
Similarly, we can verify that $W=V^-$.

If there is some vertex $w^*\in V^-$ such that $|N_\mathbf{G}(w^*)\cap V^+|<\kappa-n_2+1$,
then $(V^-\backslash \{w^*\}) \cup (N_\mathbf{G}(w^*)\cap V^+)$ is a vertex cut of $\mathbf{G}$.
Whereas $|(V^-\backslash \{w^*\}) \cup (N_\mathbf{G}(w^*)\cap V^+)|<\kappa$,
this contradiction shows that for any $w\in V^-$, $|N_\mathbf{G}(w)\cap V^+|\geq \kappa-n_2+1$.
Similarly, we can verify that for any $u^*\in V^+$, $|N_\mathbf{G}(u^*)\cap V^-|\geq \kappa-n_1+1$.

If there are two vertices $w_1$ and $w_2$ of $V^-$ such that
$N_\mathbf{G}(w_1)\cap V^+=N_\mathbf{G}(w_2)\cap V^+$ and $|N_\mathbf{G}(w_1)\cap V^+|=\kappa-n_2+1$,
then $(V^-\backslash \{w_1\cup w_2\}) \cup (N_\mathbf{G}(w_1)\cap V^+)$ is a vertex cut of $\mathbf{G}$.
Whereas $|(V^-\backslash \{w_1\cup w_2\}) \cup (N_\mathbf{G}(w_1)\cap V^+)|<\kappa$,
this contradiction shows that for any two vertices $w_3$ and $w_4$ of $V^-$ satisfying if they are all adjacent to $\kappa-n_2+1$ vertices in $V^+$, then $N_\mathbf{G}(w_3)\cap V^+\neq N_\mathbf{G}(w_4)\cap V^+$.
Similarly, we can verify that for any two vertices $u'$ and $u''$ of $V^+$ satisfying if they are all adjacent to $\kappa-n_1+1$ vertices in $V^-$, then $N_\mathbf{G}(u')\cap V^-\neq N_\mathbf{G}(u'')\cap V^-$.

We denote by $\mathbf{G}_0$ the graph obtained from $\mathbf{G}$ by connecting respectively all pairs of vertices between $V_1^+$ and $V_2^+$ and between $V_1^-$ and $V_2^-$ ($V_1^+, V_2^+, V_1^-, V_2^-$ are detailed in Claims 3.1 and 3.2).
Clearly, $\mathbf{G}_0^c$ is a bipartite graph,
and so $\lambda_n(\mathbf{G}_0^c)=-\lambda_1(\mathbf{G}_0^c)$,
where $\lambda_1(\mathbf{G}_0^c)$ is the spectral radius of $\mathbf{G}_0^c$.
Let $\Delta(\mathbf{G}_0^c)$ be the maximum degree of $\mathbf{G}_0^c$.
From the above argument, we know that $\Delta(\mathbf{G}_0^c)\leq n-\kappa-1$.
By the lemma 2.2, we know that $\lambda_1(\mathbf{G}_0^c)\leq \Delta(\mathbf{G}_0^c)$,
and so $\lambda_n(\mathbf{G}_0^c)\geq -\Delta(\mathbf{G}_0^c)\geq \kappa+1-n$.

We can easily observe that $\mathbf{B}_3^c(n_1,n_2;\kappa)$ is composed of a ($n-\kappa-1$)-regular bipartite graph and $n_1-n_2$ isolate vertices,
and so $\lambda_n(\mathbf{B}_3^c(n_1,n_2;\kappa))=\kappa+1-n$.
Therefore, $\lambda_n(\mathbf{G}_0^c)\geq \lambda_n(\mathbf{B}_3^c(n_1,n_2;\kappa))$.

From the equation $(1)$,
$x^TA(\mathbf{G})x=\sum_{v_iv_j\in E(\mathbf{G})}2x_ix_j \leq \sum_{v_iv_j\in E(\mathbf{G}_0)}2x_ix_j=x^TA(\mathbf{G}_0)x$.
As in the proof of Claim 3.1,
we can verify $\lambda_n(\mathbf{G}^c)\geq \lambda_n(\mathbf{G}_0^c)$.
Therefore, we have $\lambda_n(\mathbf{G}^c)\geq  \lambda_n(\mathbf{B}_3^c(n_1,n_2;\kappa))$.
$\Box$\vskip 3mm

{\noindent \bf Lemma 3.4.} {\it  $\lambda _n(\mathbf{B}_1^c(n_1,n_2,\kappa))> \lambda _n(\mathbf{B}_1^c(\lceil \frac{n}{2} \rceil,\lfloor \frac{n}{2} \rfloor,\kappa))$.}

{\noindent \bf Proof.}
Assume without loss of generality that $n_1 > n_2+1$.
Suppose that $U^*\subseteq V(K_{n_1})$ is a vertex cut of $\mathbf{B}_1(n_1,n_2,\kappa)$,
and that $W^*\subseteq V(K_{n_2})$ is a set composed of those vertices which are adjacent to the vertices in $U^*$.
Let $x=(x_1,x_2,\cdots,x_n)^T$ be a unit eigenvector of $\mathbf{B}_1^c(n_1,n_2,\kappa)$
with respect to $\lambda _n(\mathbf{B}_1^c(n_1,n_2,\kappa))$.
By the symmetry of $\mathbf{B}_1^c(n_1,n_2,\kappa)$,
all the vertices in $V(K_{n_1})\setminus U$ correspond to the same value $x_1$,
all the vertices in $U^*$ correspond to the same value $x_2$,
all the vertices in $V(K_{n_2})\setminus W^*$ correspond to the same value $x_3$,
and all the vertices in $W^*$ correspond to the same value $x_4$.
From the equation $(2)$, we have
$$ \left\{ \begin{array}{l}
	\lambda_nx_1=(n_2-\kappa)x_3+\kappa x_4,\\
	\lambda_nx_2=(n_2-\kappa)x_3+(\kappa-1)x_4,\\
	\lambda_nx_3=(n_1-\kappa)x_1+\kappa x_2,\\
	\lambda_nx_4=(n_1-\kappa)x_1+(\kappa-1)x_2.\\
\end{array} \right.$$

Transform the above equations into a matrix equation $(A_{n_1,n_2}-\lambda_nI_4)\widetilde{x}=0$,
where $\widetilde{x}=(x_1,x_2,x_3,x_4)^T$ and
$$A_{n_1,n_2}= \left( \begin{matrix}
	0&	  0&    n_2-\kappa&    \kappa\\
	0&    0&    n_2-\kappa&    \kappa-1\\
	n_1-\kappa&    \kappa&    0&    0\\
	n_1-\kappa&    \kappa-1&    0&    0\\
\end{matrix}  \right).$$
Let $g_{n_1,n_2}(\lambda)=\det (A_{n_1,n_2}-\lambda I_4)$.
We can compute out
\begin{equation}
	\begin{split}
		g _{n_1,n_2} ( \lambda  ) =\lambda ^4 +(2\kappa -n_1n_2-1)\lambda^2+\kappa ^2-(n_1+n_2)\kappa +n_1n_2.
		\nonumber
	\end{split}
\end{equation}
Therefore, $g _{n_1,n_2} ( \lambda  )-g _{n_1-1,n_2+1} ( \lambda  )=(n_1-n_2-1)(\lambda^2-1).$

Note that $\mathbf{B}_1^c(n_1,n_2,\kappa)$ is a bipartite graph.
It is well known that $\lambda_n(\mathbf{B}_1^c(n_1,n_2,\kappa))=-\lambda_1(\mathbf{B}_1^c(n_1,n_2,\kappa))$.
Recall, $n_1 > n_2+1$.
By Lemma 2.7, $\Delta(\mathbf{B}_1^c(n_1,n_2,\kappa))\geq n_1-1>1$,
and so by Lemma 2.2, $\lambda_n(\mathbf{B}_1^c(n_1,n_2,\kappa))<-1$.
This implies that $g_{n_1-1,n_2+1}(\lambda_n(\mathbf{B}_1^c(n_1,n_2,\kappa)))<0$.
We can observe that the function $g_{n_1-1,n_2+1}(\lambda)$ monotonically decreases when $\lambda < -\sqrt{\frac{n_1n_2+n_1-n_2-2\kappa}{2}}$,
and so $ \lambda_n(\mathbf{B}_1^c(n_1,n_2,\kappa)) > \lambda_n(\mathbf{B}_1^c(n_1-1,n_2+1,\kappa))$.
This shows $\lambda _n(\mathbf{B}_1^c(n_1,n_2,\kappa))> \lambda _n(\mathbf{B}_1^c(\lceil \frac{n}{2} \rceil,\lfloor \frac{n}{2} \rfloor,\kappa))$.
$\Box$ \vskip 3mm

{\noindent \bf Lemma 3.5.} {\it
When $n<2\kappa$, $\lambda _n(\mathbf{B}_2^c(n_1,n_2;\kappa))\geq \lambda _n(\mathbf{B}_2^c(\kappa,n-\kappa;\kappa))$,
and when $n\geq 2\kappa$, $\lambda _n(\mathbf{B}_2^c(n_1,n_2;\kappa))\geq \lambda _n(\mathbf{B}_2^c(n-\kappa+1,\kappa-1;\kappa))$.}

{\noindent \bf Proof.}
Let $y=(y_1,y_2,\cdots,y_n)^T$ be a unit eigenvector of $\mathbf{B}_2^c(n_1,n_2;\kappa)$
with respect to $\lambda _n(\mathbf{B}_2^c(n_1,n_2;\kappa))$.
By the symmetry of $\mathbf{B}_2^c(n_1,n_2;\kappa)$,
all the vertices in $V(K_{n_1})\backslash (V(M_1)\cup R_1)$ correspond to the same value $y_1$,
all the vertices in $R_1$ correspond to the same value $y_2$,
all the vertices in $V(K_{n_1})\cap M_1$ correspond to the same value $y_3$,
and all the vertices in $V(K_{n_2})$ correspond to the same value $y_4$.
From the equation $(2)$, we have
$$ \left\{ \begin{array}{l}
	\lambda_ny_1=n_2 y_4,\\
	\lambda_ny_2=0,\\
	\lambda_ny_3=(n_2-1) y_4,\\
	\lambda_ny_4=(n_1-\kappa)y_1+(n_2-1)y_3.\\
\end{array} \right.$$

Transform the above equations into a matrix equation $(A_{n_1,n_2}-\lambda_nI_4)\widetilde{y}=0$,
where $\widetilde{y}=(y_1,y_2,y_3,y_4)^T$ and
$$A_{n_1,n_2}= \left( \begin{matrix}
	0&	  0&    0&    n_2\\
	0&    0&    0&    0\\
	0&    0&    0&    n_2-1\\
	n_1-\kappa&   0&    n_2-\kappa&    0\\
\end{matrix}  \right).$$
Let $f_{n_1,n_2}(\lambda)=\det (A_{n_1,n_2}-\lambda I_4)$.
We can compute out
\begin{equation}
	\begin{split}
		f _{n_1,n_2} ( \lambda  ) =\lambda^2(\lambda^2-(n_2-1)^2-(n_1-\kappa)n_2),
		\nonumber
	\end{split}
\end{equation}
from which we obtain $\lambda _n(\mathbf{B}_2^c(n_1,n_2;\kappa))=-\sqrt{(n_2-1)^2+(n_1-\kappa)n_2}$.
Then we have $\lambda _n(\mathbf{B}_2^c(n_1-1,n_2+1;\kappa))=-\sqrt{n_2^2+(n_1-\kappa-1)(n_2)+1}$.
Recall $\kappa\leq n-2$.
By a simple computation we can determine that $\lambda _n(\mathbf{B}_2^c(n_1,n_2;\kappa))\geq \lambda _n(\mathbf{B}_2^c(n_1-1,n_2+1;\kappa))$.
This implies that
$\lambda _n(\mathbf{B}_2^c(n_1,n_2;\kappa))\geq \lambda _n(\mathbf{B}_2^c(\kappa,n-\kappa;\kappa))=\kappa +1 -n$ if $n<2\kappa$,
and $\lambda _n(\mathbf{B}_2^c(n_1,n_2;\kappa))\geq \lambda _n(\mathbf{B}_2^c(n-\kappa+1,\kappa-1;\kappa))=-\sqrt{(\kappa-2)^2+(n-2\kappa+1)(\kappa-1)}$ if $n\geq 2\kappa$.
$\Box$ \vskip 2mm


It is easy to compute that when $n< 2\kappa$,
$ \lambda _n(\mathbf{B}_3^c(n_1,n_2;\kappa))=\lambda _n(\mathbf{B}_2^c(\kappa,n-\kappa;\kappa))=\kappa+1-n$.

From Lemmas 3.2, 3.3 and 3.5,
we can easily see that the following result is true.\vskip 2mm

{\noindent \bf Theorem 3.1.} {\it When $n< 2\kappa$,
	$\mathbf{G}$ is isomorphic to $\mathbf{B}_2^c(\kappa, n-\kappa; \kappa)$ or $\mathbf{B}_3^c(n_1, n_2; \kappa)$.} \vskip 3mm

{\noindent \bf Lemma 3.6.} {\it When $n\geq 2\kappa$,
	 $ \lambda _n(\mathbf{B}_1^c(\lceil \frac{n}{2} \rceil,\lfloor \frac{n}{2} \rfloor,\kappa))<\lambda _n(\mathbf{B}_2^c(n-\kappa+1,\kappa-1,\kappa))$.}

{\noindent \bf Proof.}
From Lemma 3.4, we know $$g _{\lceil \frac{n}{2} \rceil,\lfloor \frac{n}{2} \rfloor} ( \lambda  ) =\lambda ^4 +(2\kappa -\lceil \frac{n}{2} \rceil \lfloor \frac{n}{2} \rfloor-1)\lambda^2+\kappa ^2-n\kappa +\lceil \frac{n}{2} \rceil \lfloor \frac{n}{2} \rfloor.$$
From Lemma 3.5, we know $$f _{n-\kappa+1,\kappa-1} ( \lambda  ) =\lambda^2(\lambda^2-(\kappa-2)^2-(n-2\kappa+1)(\kappa-1)).$$

Set $\phi(\lambda)=g _{\lceil \frac{n}{2} \rceil,\lfloor \frac{n}{2} \rfloor} ( \lambda  )-f _{n-\kappa+1,\kappa-1} (\lambda)$.
Then
\begin{equation}
	\begin{split}
		\phi(\lambda)=(-(\kappa ^2-n\kappa+\lceil \frac{n}{2} \rceil \lfloor \frac{n}{2} \rfloor)+\kappa-n+2)\lambda^2+\kappa ^2-n\kappa+\lceil \frac{n}{2} \rceil \lfloor \frac{n}{2} \rfloor.
		\nonumber
	\end{split}
\end{equation}
Then we can compute out the minimum root of $\phi(\lambda)$ is $\lambda_0=-\sqrt{\frac{-(\kappa ^2-n\kappa+\lceil \frac{n}{2} \rceil \lfloor \frac{n}{2} \rfloor)}{-(\kappa ^2-n\kappa+\lceil \frac{n}{2} \rceil \lfloor \frac{n}{2} \rfloor)+\kappa-n+2}}$.
Clearly, $\lambda _n(\mathbf{B}_2^c(n-\kappa+1,\kappa-1,\kappa))=-\sqrt{(\kappa-2)^2+(n-2\kappa+1)(\kappa-1)}< \lambda_0$,
and so $g _{\lceil \frac{n}{2} \rceil,\lfloor \frac{n}{2} \rfloor} (\lambda _n({B}_2^c(n-\kappa+1,\kappa-1,\kappa-1)))-f _{n-\kappa+1,\kappa-1} (\lambda _n(\mathbf{B}_2^c(n-\kappa+1,\kappa-1,\kappa)))<0$.
Thus, $g _{\lceil \frac{n}{2} \rceil,\lfloor \frac{n}{2} \rfloor} (\lambda _n(\mathbf{B}_2^c(n-\kappa+1,\kappa-1,\kappa)))<0$.
It is easy to obverse that the function $g _{\lceil \frac{n}{2} \rceil,\lfloor \frac{n}{2} \rfloor}(\lambda)$ monotonically decrease when $\lambda<\lambda _n(\mathbf{B}_1^c(\lceil \frac{n}{2} \rceil,\lfloor \frac{n}{2} \rfloor,\kappa))$,
and so $\lambda _n(\mathbf{B}_1^c(\lceil \frac{n}{2} \rceil,\lfloor \frac{n}{2} \rfloor,\kappa))<\lambda _n(\mathbf{B}_2^c(n-\kappa+1,\kappa-1,\kappa))$.   $\Box$\vskip 3mm

From Lemmas 3.1, 3.2, 3.4, 3.5 and 3.6 we can easily see that the following result is true.\vskip 2mm

{\noindent \bf Theorem 3.2.}  {\it When $n\geq 2\kappa$, $\mathbf{G}$ is isomorphic to $\mathbf{B}_1^c(\lceil \frac{n}{2} \rceil,\lfloor \frac{n}{2} \rfloor,\kappa)$.}

\end{document}